\documentclass{amsart}

\usepackage{amssymb,amsfonts, latexsym,amsmath,hhline,array,longtable}
\usepackage{pdfsync,color, comment,colortbl, mathrsfs,stmaryrd,cite,graphicx,yfonts}

\usepackage{ifpdf}
\ifpdf \usepackage[colorlinks=true, citecolor=blue, linkcolor=blue, urlcolor=blue]{hyperref} \fi

\newcommand{\cal}{\mathcal}

\newtheorem{formula}{}[section]
\newtheorem{definition}[formula]{Definition}
\newtheorem{corollary}[formula]{Corollary}
\newtheorem{remark}[formula]{Remark}
\newtheorem{lemma}[formula]{Lemma}
\newtheorem{theorem}[formula]{Theorem}
\newtheorem*{claim}{Claim}

\def\thrm{\begin{theorem}}
\def\thrml#1{\begin{theorem}\label{#1}}
\def\ethrm{\end{theorem}}
\def\rmrk{\begin{remark}}
\def\rmrkl#1{\begin{remark}\label{#1}}
\def\ermrk{\end{remark}}
\def\dfntn{\begin{definition}}
\def\dfntnl#1{\begin{definition}\label{#1}}
\def\edfntn{\end{definition}}
\def\nmrt{\begin{enumerate}}
\def\enmrt{\end{enumerate}}
\def\tm#1{\item[{\rm (#1)}]}
\def\qtnl#1{\begin{equation}\label{#1}}
\def\eqtn{\end{equation}}
\def\lmm{\begin{lemma}}
\def\lmml#1{\begin{lemma}\label{#1}}
\def\elmm{\end{lemma}}
\def\crllr{\begin{corollary}}
\def\crllrl#1{\begin{corollary}\label{#1}}
\def\ecrllr{\end{corollary}}
\def\css{\begin{cases}}
\def\ecss{\end{cases}}
\def\prf{\begin{proof}}
\def\eprf{\end{proof}}
\def\clm{\begin{claim}}
\def\eclm{\end{claim}}

\def\cA{{\cal A}}

\def\cS{{\cal S}}
\def\cT{{\cal T}}

\def\cX{{\cal X}}
\def\cY{{\cal Y}}

\def\mF{{\mathbb F}}

\def\mZ{{\mathbb Z}}

\DeclareMathOperator{\aut}{Aut}

\DeclareMathOperator{\AGL}{AGL}
\DeclareMathOperator{\AGaL}{A{\rm \Gamma}L}

\DeclareMathOperator{\fix}{fix}
\DeclareMathOperator{\GaL}{{\rm \Gamma}L}

\DeclareMathOperator{\GL}{GL}
\DeclareMathOperator{\GP}{GP}

\DeclareMathOperator{\LS}{LS}

\DeclareMathOperator{\NL}{NL}

\DeclareMathOperator{\ord}{ord}

\DeclareMathOperator{\Span}{Span}

\DeclareMathOperator{\sym}{Sym}
\DeclareMathOperator{\VLS}{VLS}
\DeclareMathOperator{\WL}{WL}

\def\und#1{{\underline{#1}}}

\begin{document}

\title[On  some generalized Paley Graphs]{The automorphism groups and identification of some Generalized Paley Graphs}
\author{}
\address{}
\author{Ilia Ponomarenko}
\address{School of Mathematical Sciences, Hebei Key Laboratory of Computational Mathematics
and Applications, Hebei Normal University, Shijiazhuang 050024, P. R. China}
\address{Steklov Institute of Mathematics at St. Petersburg, Russia}
\email{inp@pdmi.ras.ru}
\thanks{The author was supported by the grant of The Natural Science Foundation of Hebei Province (project No. A2023205045).}
\date{}

\begin{abstract}
The family of generalized Paley graphs of prime power order $q$ and degree $(q-1)/k$ is studied. It is shown that the automorphism group of a graph in this family is a subgroup of $\AGaL(1,q)$ whenever $q$ is sufficiently large relative to~$k$. Furthermore, under the same conditions, the Weisfeiler--Leman dimension of these graphs is proved to be at most~$5$. In particular, the same bound holds for the Van Lint--Schrijver graphs.	
\end{abstract}

\maketitle

\section{Introduction}
The problem of computing the automorphism group of a finite graph is known to be equivalent to the graph isomorphism problem. Nevertheless, for a particular family of graphs, determining their automorphism groups is not always directly related to the problem of identifying these graphs in the class of all graphs. Let us illustrate these two problems on the example of  the Paley graphs. 

For each prime power $q \equiv 1 \pmod{4}$, there is a unique Paley  graph, up to isomorphism. The vertices of this graph are the elements of the field $\mF_q$ with~$q$ elements, and two distinct vertices are adjacent if their difference is a square. R.~McConnell proved in 1963 that the automorphism group of this graph is a subgroup of index~$2$ in $\AGaL(1,q)$. Other proofs and a historical overview of this result can be found in the excellent survey~\cite{Jones2020}. 

On the other hand, identifying Paley graphs requires finding invariants that distinguish them from other strongly regular graphs with the same parameters; note that for a fixed $q$, there can be exponentially many of these. In \cite{Ponomarenko2020a}, it was shown that apart from finitely many cases the invariants can be taken as the triple intersection numbers of the corresponding Paley association scheme. 

In the present paper, an asymptotic solution to both problems is provided for the generalized Paley graphs defined below.

Let $k\ge 2$ be a divisor of $q-1$ such that $(q-1)/k$ is even if $q$ is odd. In \cite{LimP2009}, the authors defined the \emph{generalized Paley graph} $X=\GP(q,\frac{q-1}{k})$ with vertex set~$\mF_q$ and edges $\{x, y\}$ such that $x-y$ belongs to the subgroup of index $k$ of the multiplicative group of the field~$\mF_q$. They also provided criteria (in terms of the parameters $q$ and~$k$) for the graph $X$ to be connected and for it to be isomorphic to a Hamming graph. They also proved that if $k$ divides $p-1$, then 
\qtnl{010825a}
\aut(X)\le \AGaL(1,q). 
\eqtn
Later, in \cite{Pearce2019}, Cartesian decomposable generalized Paley graphs were completely classified.

The first of our main results extends the class of generalized Paley graphs $X$ for which inclusion~\eqref{010825a} holds true. Unlike the result in \cite{LimP2009}, we do not require that $k$ divides $p-1$, but instead impose the condition that $q$ be sufficiently large relative to~$k$;
as for an explicit lower bound on~$q$ depending on~$k$, see  Remark~\ref{240725r}.  Examples of graphs satisfying this condition can be found among the Van Lint--Schrijver graphs, as discussed later in the introduction.

\thrml{240725g}
Let $k\ge 2$ be an integer. For all $q$ sufficiently large with respect to~$k$, the automorphism group of a generalized Paley graph $X=\GP(q,\frac{q-1}{k})$ satisfies inclusion~\eqref{010825a}.
\ethrm

A natural way to identify a graph employs the Weisfeiler--Leman (WL) method, which is based on an algorithm of the same name (see, for example, \cite{Grohe2017}). In the $m$-dimensional version of this algorithm ($m$-dim WL), a canonical partition of the $m$th Cartesian power of the vertex set is constructed for a given graph $X$, $m \ge 1$. This partition defines an $(m+1)$-dimensional tensor of intersection numbers; for the case of $m=2$, these are the standard intersection numbers of the coherent configuration associated with~$X$ (see Subsection~\ref{211125a}). 

The $m$-dim WL algorithm does not \emph{distinguish} two graphs, $X$ and $X'$, if their corresponding tensors are identical. The algorithm is said to \emph{identify} a graph~$X$ if it distinguishes $X$ from all nonisomorphic graphs. The \emph{Weisfeiler--Leman dimension} of a graph $X$ is the smallest integer $m$ for which the $m$-dim WL algorithm identifies~$X$. Consequently, any graph of the Weisfeiler--Leman dimension~$m$ is identified up to isomorphism by the $(m+1)$-dimensional tensor of intersection numbers associated with the partition constructed by the $m$-dim WL algorithm.

\thrml{240725h}
Let $k\ge 2$ be an integer. The Weisfeiler--Leman dimension of a generalized Paley graph $\GP(q,\frac{q-1}{k})$ is at most $5$ for all sufficiently large $q$.
\ethrm

The proof of both Theorem~\ref{240725g} and Theorem~\ref{240725h} is based on the Babai--G\'{a}l--Wigderson theorem \cite[Lemma~3.6]{BabaiGW1999}, which provides a Weil-type bound on the number of solutions to a system of special equations over a finite field. In Section~\ref{300725a}, we slightly generalize this result to also cover inequalities of the same type. This enables us to estimate the intersection sizes of a family of vertex neighborhoods and their complements in a generalized Paley graph .

Using this technique, which is developed in Section~\ref{300725c}, we show that if $q$ is sufficiently large and $\Delta$ is the neighborhood of zero in the graph $X = \GP(q, \frac{q-1}{k})$, then:
\nmrt
\tm{a} no non-identity automorphism of the induced subgraph $X_\Delta$ leaves more than half of its vertices fixed;
\tm{b} any two distinct nonzero vertices of $X$ that do not belong to $\Delta$ cannot have the same neighbors in $\Delta$.
\enmrt
Note that the graph $X_\Delta$ is circulant:  the  cyclic group of $C_m$ with $m = \frac{q-1}{k}$, being the subgroup of index~$k$ in the group~$(\mF_q)^\times$, acts as an automorphism group of~$X_\Delta$. Then statement (a) implies that $\aut(X_\Delta)\le C_m \rtimes \aut(C_m)$; here we heavily use the theory of circulant graphs developed in \cite{EvdP2003, Evdokimov2016}. Statement~(b) then implies that the point stabilizer of $\aut(X)$ acts faithfully on $\Delta$. This is sufficient to prove Theorem~\ref{240725g}. A similar argument, translated into combinatorial language, yields the proof of Theorem~\ref{240725h}.

We conclude the introduction by considering a special family of generalized Paley graphs, namely, the Van Lint--Schrijver graphs. These graphs were first introduced and studied in~\cite{vanLintS1981}. Constructed for any two distinct primes $e \neq 2$ and $p$ such that $\ord_p(e)=p-1$, and a positive integer~$t$, they coincide with the graphs $\GP(q, \frac{q-1}{k})$ where $q = p^{(e-1)t}$ and $k = e$. We will denote them as $\VLS(p, e, t)$.

As shown in \cite[Theorem~3]{Muzychuk2021}, the Van Lint--Schrijver graphs, together with two other families (Paley graphs and Peisert graphs), are the only (up to complement) rank~$3$ graphs of prime power order~$q$ that admit an edge-transitive automorphism group contained in $\AGaL(1,q)$. More recently, it was proved in \cite[Theorem 1.2]{Skresanov2021} that, with finitely many exceptions, the automorphism group of the Van Lint--Schrijver graph of order~$q$ is contained in $\AGaL(1,q)$.\footnote{The automorphism group of the exceptional graphs was found in~\cite{GVW2024}.} It is worth noting that for all sufficiently large~$q$, this statement follows also from Theorem~\ref{240725g}. 

\crllrl{240725i}
Let $e$ and $p$ be primes such that $\ord_p(e)=p-1$. Then for all sufficiently large integers~$t$, the Weisfeiler--Leman dimension of the Van Lint--Schrijver graph $\VLS(p,e,t)$ is between $3$ and~$5$.
\ecrllr

The upper bound  in Corollary~\ref{240725i} follows immediately from Theorem~\ref{240725h}. In order to prove the lower bound, we use the fact observed in~\cite{vanLintS1981} that depending on whether $t$ is even or odd, the Van Lint--Schrijver graph $X=\VLS(p,e,t)$ is strongly regular with  Latin square  or negative Latin square parameters, see~\cite{BrouwerM2022}.  For all large prime powers~$q$, there are exponentially many of pairwise nonisomorphic strongly regular graphs of both types. The parameters of them can be chosen to be the same as those of~$X$. Therefore the $2$-dim WL algorithm does not distinguish $X$ and any of these graphs. This explains the lower bound in Corollary~\ref{240725i}. However, the question of the exact value of the Weisfeiler--Leman dimension of the Van Lint--Schrijver graphs (at least for large~$q$) remains open.

The author thanks M.~Muzychuk for helpful suggestions concerning the proof of the lower bound in Corollary~\ref{240725i}.

\section{Graphs and coherent configurations}

In this section, we introduce key definitions and basic facts concerning coherent configurations. For proofs and more detailed examples, see the monograph~\cite{CP2019}.

\subsection{Notation}
Throughout this paper, let $\Omega$ be a finite set. For any subset $\Delta \subseteq \Omega$, we denote by $1_\Delta$ the diagonal of the Cartesian square $\Delta \times \Delta$, which consists of all pairs $(\alpha, \alpha)$ where $\alpha \in \Delta$. For a binary relation $r \subseteq \Omega \times \Omega$, we define its converse as $r^* = \{(\beta, \alpha) \mid (\alpha, \beta) \in r\}$ and the neighborhood $\alpha r=\{\beta\in\Omega:\ (\alpha,\beta)\in r\}$ of the point~$\alpha$ in~$r$. Given a bijection $f$ from $\Omega$ to another set, we define the images of~$\alpha$ and~$r$ under $f$ as $\alpha^f$ and $r^f = \{(\alpha^f, \beta^f) \mid (\alpha, \beta) \in r\}$, respectively.

\subsection{Basic definitions} Let  $S$ be a partition of $\Omega\times\Omega$; we refer to 
 the elements of~$\Omega$  as points and treat the elements of~$S$  as binary relations on~$\Omega$. A pair $\cX=(\Omega,S)$ is called a \emph{coherent configuration} on $\Omega$ if the following conditions are satisfied:
\nmrt
\tm{C1}  $1_\Omega$ is a union of some relations from~$S$,
\tm{C2} for each $s\in S$, the relation $s^*$ belongs to~$S$,
\tm{C3} given any $r,s,t\in S$, the number $c_{rs}^t=|\alpha r\cap \beta s^{*}|$ does not depend on choice of the pair $(\alpha,\beta)\in t$.
\enmrt
In what follows, we write $S=S(\cX)$. The \emph{automorphism group} $\aut(\cX)$ of the coherent configuration~$\cX$ consists of all permutations $f\in\sym(\Omega)$ such that $s^f=s$ for all $s\in S$.

A set $\Delta \subseteq \Omega$ is called a \emph{fiber} of $\cX$ if the diagonal relation $1_\Delta$ is a union of some relations from $S$. Let $s \in S$ and let $\text{dom}(s) = \{\alpha \in \Omega \mid \alpha s \ne \varnothing\}$ be the domain of $s$. The size $|\alpha s|$ of the neighborhood of a point~$\alpha$ in~$s$ is constant for all $\alpha \in \text{dom}(s)$. This number is equal to the intersection number $c_{ss^*}^t$ for any $t \in S$ with $t \subseteq 1_{\text{dom}(s)}$.
 
A coherent configuration $\cX$ is said to be \emph{discrete} if each of its fibers is a singleton. It is called \emph{homogeneous} or a \emph{scheme} if $\Omega$ is its unique fiber. Given a subset $\Delta \subseteq \Omega$ that is a union of fibers, the pair $\cX_\Delta=(\Delta,S_\Delta)$, where $S_\Delta=\{s\in S \mid s\subseteq \Delta\times\Delta\}$, is a coherent configuration on $\Delta$. 
 
\subsection{Extensions and coherent closure}
The set of all coherent configurations on~$\Omega$ is equipped with a natural partial order, denoted by $\le$. A coherent configuration $\cX' = (\Omega, S')$ is said to be an \emph{extension} of $\cX = (\Omega, S)$, written as $\cX \le \cX'$, if the partition $S'$ is a refinement of the partition~$S$. That is, every relation in $S$ is a union of relations from $S'$. The maximal element with respect to this order is the discrete coherent configuration, while the minimal element is the \emph{trivial} coherent configuration. The partition $S$ of the trivial configuration consists of $1_\Omega$ and its complement in $\Omega \times \Omega$ (assuming $|\Omega| > 1$). 

The \emph{coherent closure} $\WL(T)$ of a set $T$ of binary relations on a set $\Omega$ is defined as the smallest coherent configuration $\mathcal{X}$ on $\Omega$ such that every relation in $T$ is a union of relations from $S$. The operator $\WL$ is a closure operator on the set of all partitions of $\Omega \times \Omega$, with respect to the natural partial order on these partitions. This closure can be efficiently computed by the classical Weisfeiler-Leman algorithm (see, e.g., \cite{AndresHelfgott2017}).

A special case of the coherent closure $\WL(T)$ arises when $T$ is the union of the set $S(\cX)$ for some coherent configuration $\mathcal{X}$ on $\Omega$ and the singleton relations $\{(\alpha, \alpha)\}$ for all $\alpha$ in a subset $\Delta \subseteq \Omega$. In this case, we call $\WL(T)$ the \emph{$\Delta$-extension} or an \emph{$m$-point extension} of $\mathcal{X}$, where $m = |\Delta|$. If $\Delta = \{\alpha, \beta, \ldots\}$, this extension is denoted by~$\mathcal{X}_{\alpha, \beta, \ldots}$. It holds that
$$ 
\aut(\mathcal{X}_{\alpha,\beta,\ldots})=\aut(\mathcal{X})_{\alpha,\beta,\ldots}, 
$$
where the group on the right-hand side is the pointwise stabilizer of the set $\{\alpha, \beta, \ldots\}$ in $\aut(\mathcal{X})$.

\subsection{Graphs}\label{211125a}
Unless stated otherwise, a \emph{graph} $X$ is a finite, simple, and undirected graph. Its \emph{order} is the number of vertices. We write $\alpha \sim \beta$ to indicate that two vertices~$\alpha$ and~$\beta$ are adjacent in~$X$. The subgraph of $X$ induced by a vertex subset~$\Delta$ is denoted by $X_{\Delta}$. The automorphism group of $X$ is denoted by $\aut(X)$.

The \emph{coherent configuration} $\WL(X)$ associated with $X$ is the coherent closure $\WL(T)$ where $T$ is the singleton consisting the edge set of $X$. A fundamental property we will use further is that $\aut(X)=\aut(\WL(X))$.

\section{Normal circulant graphs}
A graph $X$ is called \emph{circulant} if its vertex set can be identified with the elements of a cyclic group $G$  acting as a regular subgroup of $\aut(X)$ via right multiplication; in other words, $X$ is isomorphic to a Cayley graph of a cyclic group.  Consequently, $\aut(X)$ is a transitive permutation group, and the coherent configuration $\WL(X)$ associated with~$X$ is a scheme.

A circulant graph $X$ is said to be \emph{normal} if its automorphism group $\aut(X)$ is a subgroup of the semidirect product $G \rtimes \aut(G)$. An equivalent condition is that the stabilizer of the identity element of $G$ in $\aut(X)$ 
is a subgroup of $\aut(G)$. For instance, the complete graph $K_n$ is a normal circulant if and only if $n \le 3$.

A complete characterization of normal circulant graphs was obtained in \cite{EvdP2003}. The following statement is a straightforward consequence of that characterization.

\lmml{200725b}
The scheme of  a normal circulant graph has a discrete $2$-point extension.
\elmm
\prf
Let $\cX$ be the scheme associated with some normal circulant graph, and let~$\alpha$ be an arbitrary point of $\cX$. By \cite[Corollary 6.2]{EvdP2003}, the coherent configuration $\mathcal{Y} = \mathcal{X}_{\alpha}$ has a \emph{regular point}~$\beta$. By definition, this means that $|\beta s| \le 1$ for all relations $s \in S(\mathcal{Y})$. Since the set of all fibers of the coherent configuration~$\mathcal{Y}_\beta$ is equal to~$\{\beta s:\ s \in S(\mathcal{Y})\}$ (see \cite[Lemma 3.3.5]{CP2019}), this implies that the $2$-point extension $\mathcal{X}_{\alpha,\beta} = \mathcal{Y}_\beta$ is discrete.
\eprf

The proof of the lemma below uses the theory of Schur rings (S-rings) over a finite group, as developed in \cite{EvdP2003, Evdokimov2016}. This theory provides a useful algebraic framework for studying schemes that admit a regular automorphism group.

Let $\cX$ be a \emph{Cayley scheme} on a group $G$. This means that the (regular) permutation group induced by right multiplication of~$G$ is a subgroup of the group~$\aut(\cX)$. With each Cayley scheme on~$G$ one can associate an S-ring $\cA$ over $G$. It is a subring of the group ring $\mZ G$,  such that
$$
\cA=\Span_\mZ\{\und{X}:\ X\in\cS\}
$$
where $\und{X}=\sum_{x\in X}x$ and  $\cS=\{\alpha s:\ s\in S(\cX)\}$. The partition $\cS$ of the group~$G$ is denoted by~$\cS(\cA)$ and the classes of it are called the \emph{basic sets} of~$\cA$.

The correspondence $\cX\mapsto\cA$ between the Cayley schemes $\cX$ on~$G$ and the S-rings~$\cA$ over $G$ respects the inclusion relation and preserves the automorphism groups, i.e., $\aut(\cX)=\aut(\cA)$.  A full treatment of Schur ring theory is beyond the scope of this paper, so we refer the reader to the seminal papers mentioned above.

\lmml{190725a}
Any non-normal circulant graph of order $n$ has an automorphism $\sigma$ such that 
\qtnl{311025b}
n/2\le\fix(\sigma)\le 2n/3, 
\eqtn
where $\fix(\sigma)$ is the number of vertices fixed by~$\sigma$.
\elmm
\prf
Let $X$ be a circulant graph of order $n$, and let $\cX=\WL(X)$ be the scheme associated with~$X$. The vertices of~$X$ are identified with the elements of a cyclic group $G$ of order~$n$ and hence $\cX$ is a Cayley scheme on~$G$. Let $\cA$ be the S-ring over~$G$ corresponding to~$\cX$. Then 
\qtnl{311025a}
\aut(X)=\aut(\cX)=\aut(\cA).
\eqtn

Suppose first that the S-ring $\cA$ has a nontrivial radical in the sense of~\cite{Evdokimov2016}. By Theorem~5.2 of that paper, this implies that $\cA$ is the $U/L$-wreath product for some subgroups $U$ and $L$ of $G$ such that $1 < L \le U < G$. This structural property means that $U$ and $L$ are unions of basic sets of~$\cA$ and any basic set of~$\cA$ that is not a subset of $U$ must be a union of $L$-cosets.

Let us define a new partition $\cS'$ of $G$ consisting of singletons for elements within~$U$ and $L$-cosets for elements outside $U$. This partition $\cS'$ is a refinement of the partition $\cS(\cA)$ of~$G$. Furthermore, $\cS'$ coincides with the set of basic sets of an S-ring over~$G$, which we denote by $\cA' = \mZ(G,U/L)$, as defined in Subsection~8.2 of~\cite{Evdokimov2016}. It follows that $\cA' \supseteq \cA$. This inclusion implies that any automorphism of $\cA'$ is also an automorphism of $\cA$, i.e., 
\qtnl{210725f}
\aut(\cA')\le\aut(\cA).
\eqtn 

By Theorem~8.2 of the same paper, the automorphism group $\aut(\cA')$ is precisely the set of permutations of the form $\sigma_{t,g} \in \sym(G)$, where $g \in G$ is an arbitrary element and $t: G/U \to L$ is an arbitrary function. These permutations act on $G$ as follows:
$$
x^{\sigma_{t,g}} = x \cdot t(Ux) \cdot g, \qquad \text{for all } x \in G.
$$
In view of \eqref{210725f}, any such permutation $\sigma_{t,g}$ is an automorphism of $\cA$.

Now, let us consider the number of fixed points for these automorphisms. Let $m = |G/U| = |G|/\,|U|$. We can choose $g=1$ and define a function $t: G/U \to L$ that maps exactly $\lceil m/2 \rceil$ cosets to the identity element $1 \in L$. The fixed points of the permutation $\sigma=\sigma_{t,1}$ are the elements $x \in G$ for which $x^{\sigma_{t,1}} = x$, which is equivalent to $x \cdot t(Ux) = x$, or simply $t(Ux)=1$.

We analyze the number $\fix(\sigma)$ in two cases based on the parity of $m$.
\begin{enumerate}
\item If $m$ is even, say $m=2k$, then $\lceil m/2 \rceil = k$. The number of cosets for which $t(Ux)=1$ is exactly $k$. Since each coset contains $|U|$ elements, the total number of fixed points is $\fix(\sigma)=k|U| = \frac{m}{2}|U| = \frac{|G|}{2}$.
\item If $m$ is odd, say $m=2k+1$ for some $k \ge 1$, then $\lceil m/2 \rceil = k+1$. The number of cosets for which $t(Ux)=1$ is $k+1$. It follows that $\fix(\sigma) =(k+1)|U| = \frac{m+1}{2}|U|$. The ratio of $\fix(\sigma)$  to the total number of points is $\frac{(k+1)|U|}{|G|} = \frac{k+1}{m} = \frac{k+1}{2k+1}$. This ratio lies in the interval $[\frac{1}{2}, \frac{2}{3}]$ for any $k \ge 1$.
\end{enumerate}
In both cases, the inequalities \eqref{311025b} hold.

Now, let us consider the case where the radical of the S-ring $\cA$ is trivial. According to Theorem~6.4 of~\cite{EvdP2003}, any such S-ring is the tensor product of a normal S-ring and trivial S-rings.\footnote{An S-ring over $G$ is said to be trivial if the basic sets of it are $\{1\}$ and $G\setminus\{1\}$ (if $|G|>1$.)} Since a trivial S-ring is normal if and only if its degree is at most $3$, and the tensor product of normal S-rings is also normal, it follows that one of two cases must hold:
\begin{enumerate}
	\item the S-ring $\cA$ itself is normal,
	\item $\cA = \cT_m \otimes \cA'$, where $\cT_m$ is a trivial S-ring of degree $m \ge 4$ and $\cA'$ is another S-ring.
\end{enumerate}
Assuming the graph $X$ is non-normal, we can exclude the first case, see~\eqref{311025a}. In the second case, the automorphism group of $\cA$ is the direct product of the automorphism groups of its factors:
$$
\aut(\cA) = \aut(\cT_m) \times \aut(\cA') \ge \sym(m) \times \{1\},
$$
where $\sym(m)$ is the full symmetric group on $m$ points, and $\{1\}$ represents the identity subgroup acting on the domain of $\cA'$. This is a subgroup of the automorphism group $\aut(\cA) = \aut(X)$.

To calculate the number of fixed points, let $\sigma_1$ be any permutation in $\sym(m)$  that leaves exactly $\lceil m/2 \rceil$ points fixed, i.e., $\fix(\sigma_1)=\lceil m/2 \rceil$.  We can construct a permutation $\sigma = (\sigma_1, 1) $ belonging to~$\sym(m) \times \{1\}$. 
In particular,  $\sigma$ belongs to the group $\aut(\cA) = \aut(X)$. 

If $m$ is even, then the number of fixed points of $\sigma_1$ is $m/2$. Thus, 
$$
\fix(\sigma)=\frac{m}{2} \cdot \frac{n}{m} = \frac{n}{2}.
$$
If $m=2k+1$ for some $k\ge 1$, then the number of fixed points of $\sigma_1$ is $k+1$. Thus,
$$
\fix(\sigma)=(k+1) \cdot \frac{n}{m} = \frac{k+1}{2k+1}n,
$$
which completes the proof, as the number of fixed points falls within the required range.
\eprf

\section{The Babai--G\'{a}l--Wigderson theorem}\label{300725a}

In this section, we slightly generalize a result deduced in~\cite[Lemma~3.6]{BabaiGW1999} from the Weil theorem on the number $\mF_q$-rational points of absolutely irreducible curves defined over the field~$\mF_q$. As noted therein, it ``provides a tight control on the intersection sizes of vertex neighborhoods in the Paley-type graphs''. Note that the latter graphs are closely related with the generalized Paley graphs. We quote the mentioned result below.

\thrml{160725a}
Let $\alpha_1,\ldots, \alpha_t$  be distinct elements of the finite field~$\mF_q$ and $k$ be a divisor of $q-1$ $(k,t\ge 2)$.  Then  the  number of solutions to the system of equations 
\qtnl{170725s}
(\alpha_i+x)^{\frac{q-1}{k}}=1\qquad(i=1,\ldots,t)
\eqtn
 is between $\frac{q}{k^t}-t\sqrt{q}$ and $\frac{q}{k^t}+t\sqrt{q}$. 
\ethrm

The proof of our main theorems in Section~\ref{030825a} relies on an estimate for the number of solutions to a specific system of equations and inequalities. We consider a system composed of $t$ equations of the form given in~\eqref{170725s} and $n-t$ inequalities of the form $(\alpha_i+x)^{\frac{q-1}{k}}\ne 1$. The full statement is provided in Lemma~\ref{160725a1}. We point out that while this estimate is sufficient for our purposes, the leading constant in the error term is likely not sharp.

\lmml{160725a1}
Let $\alpha_1,\ldots, \alpha_n$  be distinct elements of the finite field~$\mF_q$ and $k$ be a divisor of $q-1$ $(n, k\ge 2)$.  For $t \in \{2, \ldots, n\}$, let $N$ be the number of solutions $x \in \mF_q$ to the system defined by the equations~\eqref{170725s}  and the inequalities
\qtnl{170725r}
(\alpha_i+x)^{\frac{q-1}{k}}\ne 1\qquad (i=t+1,\ldots,n).
\eqtn 
Then $N$ satisfies the bound
\qtnl{311025u}
\left| N - q \left(\frac{1}{k}\right)^t \left(1-\frac{1}{k}\right)^{n-t} \right| \le n\,2^{(n-t)^2}\sqrt{q}. 
\eqtn
\elmm

\prf
For each subset $T \subseteq \{1,\ldots,n\}$, let $N(T)$ denote the number of solutions $x \in \mF_q$ to the system of $t:=|T|$ equations given in~\eqref{170725s} for $i \in T$, and $s:=n-t$ inequalities given in~\eqref{170725r} for $i \notin T$. We will prove the following estimate by induction on $t$, counting downwards from $n$ to $2$:
$$
\left| N(T) - q \left(\frac{1}{k}\right)^t \left(1-\frac{1}{k}\right)^s \right| \le n\,2^{s^2}\sqrt{q}.
$$
For $t=n$, we have $s=0$. The set $T$ is unique, $T=\{1,\ldots,n\}$, and there are no inequalities. The required estimate follows directly from the bound in Theorem~\ref{160725a}.

Assume the statement holds for all subsets of size greater than $t$. Let us consider a set $T$ of size $t < n$, so $s \ge 1$. Without loss of generality, let $T = \{1, \ldots, t\}$. The number of solutions $N_0(T)$ to the system of equations specified by $T$ (without any inequalities)   can be expressed as follows:
$$
N_0(T) = \sum_{U \subseteq \{t+1,\ldots,n\}} N(T \cup U).
$$
Rearranging the equation above, we get
$$
N(T) = N_0(T) - \sum_{\substack{U \subseteq \{t+1,\ldots,n\} \\ U \ne \varnothing}} N(T \cup U).
$$
By Theorem~\ref{160725a}, we know that $|N_0(T)-q(1/k)^t |\le t\sqrt{q}$. For the sum, we apply the inductive hypothesis to each term $N(T \cup U)$ for $|U|=u \ge 1$. Then
\qtnl{311025t}
|N(T)-N_1q|\le N_2\sqrt{q},
\eqtn
where
$$
N_1= \left(\frac{1}{k}\right)^t-\sum_{u=1}^{s}\binom{s}{u}
\left(\frac{1}{k}\right)^{t+u}\,\left(1-\frac{1}{k}\right)^{s-u} 
$$
and 
$$
N_2= t+\sum_{u=1}^{s}\binom{s}{u}n\,2^{(s-u)^2}.
$$
We calculate the quantity $N_1$ as follows:
$$
N_1= \left(\frac{1}{k}\right)^t-\left(\frac{1}{k}\right)^t\cdot\sum_{u=1}^{s}\binom{s}{u}
\left(\frac{1}{k}\right)^{u}\,\left(1-\frac{1}{k}\right)^{s-u} =
$$
$$
\left(\frac{1}{k}\right)^t\left( 
1+  
\left(1-\frac{1}{k}\right)^{s} -\sum_{u=0}^{s}\binom{s}{u}
\left(\frac{1}{k}\right)^{u}\,\left(1-\frac{1}{k}\right)^{s-u} 
\right)=
$$
$$
\left(\frac{1}{k}\right)^t\left( 
1+  
\left(1-\frac{1}{k}\right)^{s} -1\right)=\left(\frac{1}{k}\right)^t\,\left(1-\frac{1}{k}\right)^s.
$$
The error term, however, needs to be estimated more carefully. A valid upper bound can be obtained as follows:
$$
N_2 = t + n \sum_{u=1}^{s} \binom{s}{u} 2^{(s-u)^2} \le n+ n \cdot 2^{(s-1)^2}\cdot  \sum_{u=1}^{s} \binom{s}{u}=
$$
$$
n\cdot(1+  2^{(s-1)^2}\cdot (2^s-1))\le n\cdot (1+2^{(s-1)^2+s}-  2^{(s-1)^2})\le n\cdot 2^{s^2}.
$$
The final bound $n\cdot 2^{s^2}$ is a loose but valid upper bound for the error term. Now combining~\eqref{311025t} and the obtained estimates for~$N_1$ and~$N_2$, we get the required inequality~\eqref{311025u}.
\eprf

\section{Distinguishing vertices in  generalized Paley graphs}\label{300725c}

Let $X$ be a graph. A vertex $\alpha$ \emph{distinguishes} vertices $\beta$ and $\gamma$ (in $X$) if $\alpha$ is adjacent to exactly one of $\beta$ and $\gamma$. Let us assume that $X=\GP(q,\frac{q-1}{k})$. Then, $\alpha$ distinguishes $\beta$ and $\gamma$ if and only if $\alpha$ is a solution to one of the following two systems:
$$
\begin{cases}
	(x-\beta)^{\frac{q-1}{k}}=1\\
	(x-\gamma)^{\frac{q-1}{k}}\ne 1,
\end{cases}
\text{ or }
\begin{cases}
	(x-\beta)^{\frac{q-1}{k}}\ne 1\\
	(x-\gamma)^{\frac{q-1}{k}}=1.
\end{cases}
$$
Consequently, the total number of such vertices $\alpha$ is given by the sum $N^+ + N^-$, where $N^+$ and $N^-$ are the numbers of solutions to the left- and right-hand side systems, respectively. Note that these two systems have no common solutions. Using these remarks and Lemma~\ref{160725a1}, we will now prove the main result of this section.

\lmml{230725a}
 Let $X=\GP(q,\frac{q-1}{k})$  and let $\Delta$ be the neighborhood  of the vertex~$0$. Assuming that $q$ is sufficiently large with respect to $k$, the following holds:
 \nmrt
 \tm{1} any two distinct nonzero vertices of $X$ that do not belong to $\Delta$ are distinguished by some vertex in~$\Delta$,
 \tm{2} if $S$ is a set of $k$ pairwise disjoint $2$-vertex subsets of  $\Delta$, then  more than half of the vertices in $\Delta$ distinguish two vertices of some $2$-subset in~$S$.
  \enmrt
 \elmm
 
 \rmrkl{240725r}
 The proof shows that statements {\rm (1)}  and {\rm (2)} holds if $\sqrt{q}$ is larger than $ \frac{5k^3}{k-1}$ and  $10(2k^2+k)2^{k^2+2k}$, respectively.
 \ermrk
 
 \prf
Note that the assumption on~$q$ ensures that $X$ is connected (see~\cite[Theorem~1.2, statement~(1)]{LimP2009}). We prove statement (1) by contradiction. Let $\beta$ and $\gamma$ be two distinct nonzero vertices of $X$ that are not in $\Delta$. Assume, to the contrary, that no vertex in $\Delta$ distinguishes them. This implies that $\beta$ and $\gamma$ have the same set of neighbors in $\Delta$. Let us denote this set by $\Delta(\beta, \gamma) \subseteq \Delta$.

Since every vertex $\alpha \in \Delta(\beta, \gamma)$ is adjacent to both $\beta$ and $\gamma$, it must satisfy the system of two equations given by $(\alpha-\beta)^{\frac{q-1}{k}} = 1$ and $(\alpha-\gamma)^{\frac{q-1}{k}} = 1$. By applying Lemma~\ref{160725a} with $t=2$, we get a lower bound on the size of this set:
$$|\Delta(\beta, \gamma)| \ge q \left(\frac{1}{k}\right)^2 - 2\sqrt{q}.$$

On the other hand, every vertex $\alpha \in \Delta(\beta, \gamma)$ is also in $\Delta$, which means it is adjacent to the vertex $0$. Therefore, each $\alpha \in \Delta(\beta, \gamma)$ must also satisfy the system of three equations given by $(\alpha-0)^{\frac{q-1}{k}} = 1$, $(\alpha-\beta)^{\frac{q-1}{k}} = 1$, and $(\alpha-\gamma)^{\frac{q-1}{k}} = 1$. Applying Lemma~\ref{160725a} again, this time with $t=3$, we obtain an upper bound:
$$|\Delta(\beta, \gamma)| \le q \left(\frac{1}{k}\right)^3 + 3\sqrt{q}.$$

Combining these two inequalities, we have
$$q \left(\frac{1}{k}\right)^2 - 2\sqrt{q} \le |\Delta(\beta, \gamma)| \le q \left(\frac{1}{k}\right)^3 + 3\sqrt{q}.$$
This leads to the inequality
$$q \left(\frac{1}{k^2} - \frac{1}{k^3}\right) \le 5\sqrt{q},$$
which simplifies to $\sqrt{q} \le \frac{5k^3}{k-1}$. This contradicts the assumption that $q$ is sufficiently large with respect to $k$. Thus, the initial assumption must be false, and statement~(1) holds.

To prove statement (2), let $S$ be a set of $k$ pairwise disjoint $2$-vertex subsets of~$\Delta$, say $S = \{\{\alpha_i, \beta_i\} : i=1,\ldots,k\}$. For each $i$, let $A_i$ be the set of all vertices in $\Delta$ that distinguish the pair $\{\alpha_i, \beta_i\}$. We aim to prove that
\qtnl{190725b}
|A_1\cup\ldots\cup A_k|\ge \frac{|\Delta|}{2}.
\eqtn
To this end, we first estimate the size of the intersections of these sets. For any nonempty subset $T = \{i_1, \ldots, i_t\} \subseteq \{1, \ldots, k\}$, let $M(T) = |A_{i_1} \cap \cdots \cap A_{i_t}|$. We will then apply the Principle of Inclusion-Exclusion.

A vertex $x$ is in the intersection $A_{i_1} \cap \cdots \cap A_{i_t}$ if it belongs to $\Delta$ and distinguishes each pair $\{\alpha_{i_j}, \beta_{i_j}\}$ for $j=1, \ldots, t$. This means $x$ is adjacent to $0$, and for each $j \in T$, it is adjacent to exactly one of $\alpha_{i_j}$ and $\beta_{i_j}$. This leads to $2^t$ mutually exclusive systems of equations and inequalities. Each system consists of $t+1$ equations (one for vertex $0$ and one for exactly one vertex in each pair) and $t$ inequalities (for the other vertex in each pair). By applying Lemma~\ref{160725a} to these systems, we can estimate $M(T)$. Specifically, for any such system, the number of solutions is approximately $q (\frac{1}{k})^{t+1} (1-\frac{1}{k})^t$. Summing over all $2^t$ systems and using the error bound, we have
\qtnl{190725d}
\frac{q}{k}M_1(t)-M_2(t)\sqrt{q}\le M(T)\le \frac{q}{k}M_1(t)+M_2(t)\sqrt{q},
\eqtn
where $M_1(t) = 2^t(\frac{1}{k})^t(1-\frac{1}{k})^t = (\frac{2(k-1)}{k^2})^t$ and $M_2(t)=2^t(2t+1)2^{t^2}=(2t+1)2^{t^2+t}$.

Now we apply  a form of the Principle of Inclusion-Exclusion, to obtain a lower bound for the size of the union:
$$
|A_1 \cup \ldots \cup A_k| \ge \sum_{r=1}^k (-1)^{r+1} \sum_{|T|=r} M(T).
$$
Using our estimates for $M(T)$, we get a lower bound:
$$|A_1 \cup \ldots \cup A_k| \ge \frac{q}{k} \sum_{r=1}^k (-1)^{r+1} \binom{k}{r} M_1(r) - \sqrt{q} \sum_{r=1}^k \binom{k}{r} M_2(r).$$
Let the first sum be $M_1'$ and the second be $M_2'$. We can calculate $M_1'$ exactly:
$$M_1' = \sum_{r=1}^k (-1)^{r+1} \binom{k}{r} \left(\frac{2(k-1)}{k^2}\right)^r = 1 - \left(1 - \frac{2(k-1)}{k^2}\right)^k.$$
Using the inequality $(1 - a/x)^x < e^{-a}$ with $x=k$ and $a = 2(k-1)/k$, we get:
$$\left(1 - \frac{2(k-1)}{k^2}\right)^k < e^{-\frac{2(k-1)}{k}} \le e^{-1} < 0.4.$$
This yields $M_1' > 1 - 0.4 = 0.6$. Now we estimate $M'_2$ very roughly from above by $\sum_{r=1}^k \binom{k}{r}(2r+1)2^{r^2+r}< (2k+1)2^{k^2+2k}$. For $q$ large enough with respect to $k$, the last term is less than   $0.1\cdot \frac{\sqrt{q}}{k}$. This yields $M_2'< 0.1\frac{\sqrt{q}}{k}$. 

Thus, if $q$ is large enough, the lower bound for the union is given by
$$
|A_1 \cup \ldots \cup A_k| > \frac{q}{k}M_1' - M_2'\sqrt{q}\ge\frac{q}{k} 0.6 -0.1\frac{\sqrt{q}}{k}\sqrt{q}\ge \frac{1}{2}|\Delta|,
$$
which proves inequality~\eqref{190725b}.
\eprf

\section{Proof of the main results}\label{030825a}

\subsection{Preliminaries} Let $X = \GP(q, \frac{q-1}{k})$ be a generalized Paley graph, where~$q$ and $k$ are as defined in Lemma~\ref{230725a}. This choice ensures that $X$ is connected and is not a Hamming graph (see~\cite[Theorem 1.2, statements (1) and (2)]{LimP2009}). Let $\mathcal{X} = \WL(X)$ be the scheme associated with~$X$, and let $\Delta$ be the neighborhood of the vertex~$0$ in~$X$. Our first goal is to derive two consequences from Lemma~\ref{230725a}.

\lmml{250725b} 
The circulant graph $X_\Delta$ is normal.
\elmm
\prf
Assume for contradiction that $X_\Delta$ is not normal. Let $n=|\Delta|$. By Lemma~\ref{190725a}, there exists an automorphism $f\in\aut(X_\Delta)$ such that the set of fixed points $\Delta_0=\{\delta\in\Delta:\delta^f=\delta\}$ satisfies $\frac{n}{2}\le |\Delta_0|\le \frac{2n}{3}$. Let $\Delta_1=\Delta\setminus\Delta_0$  be the set of points not fixed by~$f$. Then 
$$
|\Delta_1|=n-|\Delta_0|\ge n-\frac{2n}{3}=\frac{n}{3}.
$$ 
Since $n=\frac{q-1}{k}$, our assumptions on $q$ and $k$ implies that  there exist $k$ pairwise disjoint two-vertex subsets  $\{\delta,\delta^f\}\subseteq \Delta_1$ such that $\delta\ne \delta^f$. Let $S$ be the set of these subsets, $|S|=k$.

Let $\Delta^*$ be the set of all vertices in~$\Delta$ that  distinguish some   vertices $\delta$ and~$\delta^f$  such that $\{\delta,\delta^f\}\in S$. That is, for each $\alpha\in \Delta^*$, there exists a subset $\{\delta,\delta^f\}\in S$ such that $\alpha$ is adjacent to one vertex in it but not the other. By Lemma~\ref{230725a}(2), we have $|\Delta^*|>n/2$.

Since $|\Delta_0|\ge n/2$ and $|\Delta^*|>n/2$, the intersection $\Delta_0\cap \Delta^*$ must be nonempty. Let $\alpha\in \Delta^*\cap\Delta_0$.  Then $\alpha^f=\alpha$. Since $\alpha\in \Delta^*$, there exists a pair $\{\delta,\delta^f\}\in S$ that is distinguished by $\alpha$. Without loss of generality, we may assume that $\delta\sim\alpha$ and $\delta^f\not\sim\alpha$.

The vertex $\alpha$ is fixed by the automorphism $f$, so $\alpha^f=\alpha$. The edge $\{\delta,\alpha\}$ in~$X$ must be mapped to an edge. However, the image of this edge under $f$  would be $\{\delta^f,\alpha^f\}=\{\delta^f,\alpha\}$. But we assumed that $\delta^f\not\sim\alpha$, which means $\{\delta^f,\alpha\}$ is a non-edge. This is a contradiction, as an automorphism must preserve adjacency. Thus, the initial assumption that $X_\Delta$ is not normal must be false.
\eprf

\lmml{250725a} 
The $\Delta$-extension $\cY$ of  the coherent configuration $\cX_0$  is discrete.
\elmm
\prf
Assume for contradiction that the coherent configuration $\cY$ is not discrete. This means there exist two distinct vertices, $\alpha$ and $\beta$, both not adjacent to~$0$, that belong to the same fiber of $\cY$.

A key property of the $\Delta$-extension  is that for any $\delta\in\Delta$, the singleton $\{\delta\}$  is a fiber of $\cY$. This implies that for any $\delta\in\Delta$, the adjacency relations of $\alpha$ and $\beta$ with~$\delta$ must be identical. That is, $\delta\sim\alpha$ if and only if $\delta\sim\beta$.

This statement means that no vertex in $\Delta$ can distinguish $\alpha$ and $\beta$. However, this directly contradicts Lemma~\ref{230725a}(1), which provides that for any pair of distinct vertices not adjacent to~$0$, there must be a vertex in $\Delta$ that distinguishes them. Therefore, our initial assumption is false, and $\cY$ must be discrete.
\eprf

\subsection{Proof of Theorem~\ref{240725g}} 
Recall that $X$ is a connected graph not isomorphic to a Hamming graph. In accordance with \cite[Theorem~1.2, statement~(3)]{LimP2009}, we have
$$
\mF^+\rtimes M\le \aut(X)\le \AGL(d,p),
$$
where $M$ is a subgroup of index $k$ in the group $\mF^\times$, $p$ is the characteristic of the field $\mF$, and $q=p^d$. It follows that the stabilizer $\aut(X)_0$ of the vertex $0$ in $\aut(X)$ satisfies $M \le \mathrm{Aut}(X)_0 \le \GL(d,p)$. Due to the connectivity of $X$, the group $M$ is an irreducible subgroup of a Singer cycle of  $\GL(d,p)$.

We now show that $M$ is a normal subgroup of $\aut(X)_0$. This is the central step toward proving that $\aut(X)\le \AGaL(1,q)$. Assuming that
\qtnl{270725a}
M\trianglelefteq \aut(X)_0,
\eqtn
we see that $\aut(X)_0$ must be a subgroup of the normalizer of $M$ in $\GL(d,p)$. By a known result (see \cite[Theorem~II.7.3a, p.~187]{Huppert1967}), this normalizer is contained in $\AGaL(1,q)$. Thus, we have $\aut(X)_0 \le \AGaL(1,q)$, which gives us the required inclusion:
$$
\aut(X)\le \mF^+\rtimes\aut(X)_0\le  \mF^+\rtimes\GaL(1,q)=\AGaL(1,q).
$$

To prove~\eqref{270725a}, we make use of Lemmas~\ref{250725a} and~\ref{250725b}. By Lemma~\ref{250725a}, the $\Delta$-extension of the coherent configuration $\cX_0$ is discrete. Consequently, the pointwise stabilizer of $\Delta$ in $\aut(X)_0$ is trivial. This implies that the action of the group $\aut(X)_0$ on the set $\Delta$ is faithful, leading to the embedding:
\qtnl{270725r}
\aut(X)_0\cong (\aut(X)_0)^\Delta\le\aut(X_\Delta).
\eqtn
By Lemma~\ref{250725b}, the graph $X_\Delta$ is normal. A key property of normal circulant graphs is that the regular cyclic subgroup $M^\Delta$ is a normal subgroup of $\aut(X_\Delta)$. Therefore, by the isomorphism in~\eqref{270725r}, the subgroup $M$ is normal in $\aut(X)_0$. This completes the proof of~\eqref{270725a}.

\subsection{Proof of Theorem~\ref{240725h}}
By Lemma~\ref{250725b}, the graph $X_\Delta$ is a normal circulant graph. According to Lemma~\ref{200725b}, there exist vertices $\alpha, \beta \in \Delta$ such that the coherent configuration $\WL(X_\Delta)_{\alpha,\beta}$ is discrete. On the other hand, we have the following relationships between coherent configurations:
$$
(\mathcal{X}_{0,\alpha,\beta})_\Delta \ge ((\mathcal{X}_0)_\Delta)_{\alpha,\beta} \ge \WL(X_\Delta)_{\alpha,\beta}.
$$
Since $\WL(X_\Delta)_{\alpha,\beta}$ is discrete, these inclusions imply that $((\cX)_{0,\alpha,\beta})_\Delta$ must also be a discrete coherent configuration. Since $\cX_{0,\alpha,\beta}\ge \cX_0$, this implies that  the coherent configuration $\cX_{0,\alpha,\beta}$  is larger than or equal to the $\Delta$-extension of $\cX_0$. By Corollary~\ref{250725a}, the $\Delta$-extension of $\cX_0$ is discrete. Since any coherent configuration that is larger than or equal to a discrete one must also be discrete, we conclude that $\cX_{0,\alpha,\beta}$ is discrete.

According to the terminology  in~\cite{CaiGGP2025}, the discreteness of $\mathcal{X}_{0,\alpha,\beta}$ for some $\alpha, \beta \in \Delta$ implies that the base number of the coherent configuration $\mathcal{X} = \WL(X)$ is at most 3. Finally, applying \cite[Corollary~3.2]{CaiGGP2025}, we can conclude that the Weisfeiler-Leman dimension of the graph~$X$ is at most~$5$.

\subsection{Proof of Corollary~\ref{240725i}} The upper bound  follows from Theorem~\ref{240725h}. To prove the lower bound, let $X=\VLS(p,e,t)$  the Van Lint--Schrijver graph, 
$$
q=p^d,\quad d=(e-1)t,\quad  n=\sqrt{q},\quad m=\frac{q-1}{e(n-\varepsilon)},
$$
where $\varepsilon=(-1)^t$. In accordance with~\cite{vanLintS1981} and~\cite{BrouwerM2022},  $X$ is a strongly regular graph with  Latin square  parameters $\LS_m(n)$ or with negative Latin square parameters~$\NL_m(n)$, depending on whether  $t$ is even or odd.

\lmml{070825a}
There are at least $\binom{n+\varepsilon}{m}$ distinct strongly regular Cayley graphs  of the group~$\mF_q^+$, that have the same parameters as the graph~$X$. 
\elmm
\prf
Let $\cX$ be the amorphic scheme  described in \cite[Subsection~5.4]{vanDamM2010}. The point set of $\cX$ equals $\mF_q$ and the set $S=S(\cX)$ has exactly  $r=n+\varepsilon$ nondiagonal symmetric relations invariant with respect to the group~$\mF_q^+$. The union of any $m$ of these relations is the edge set of a strongly regular graph with  parameters $\LS_m(n)$ or~$\NL_m(n)$, depending on whether  $t$ is even or odd. Since there are exactly $\binom{n+\varepsilon}{m}$ ways to choose $m$ nondiagonal relations  from $S$, we are done.
\eprf

\lmml{070825b}
Among the graphs in Lemma~\ref{070825a}, there are at most $p^{d^2}$ graphs isomorphic to~$X$.
\elmm
\prf
First, we note that $A=\mF_q^+$ is the only regular elementary abelian subgroup of the group $\AGaL(1,q)$. Indeed, this is obvious if $d$ is coprime to $p$. Assume that $p$ divides~$d$ and 
$B$ is another regular elementary abelian subgroup of the group~$\AGaL(1,q)$. The Sylow $p$-subgroup of the group $\AGaL(1,q)/A=\GL(1,q)$ is cyclic. It follows that the elementary abelian group $BA/A\le \AGaL(1,q)/A$ is generated by some $\sigma\in\aut(\mF_q)$. In particular, $BA/A$ is cyclic, and hence 
$$
|A\cap B|=p^{d-1}.
$$ 
The automorphism $\sigma$ fixes all elements of $A\cap B$. Consequently, they belong to a subfield  of~$\mF_q$. However, the order of this subfield is less than or equal to $p^{d/2}$. Thus $d-1\le d/2$, which is impossible for $d>2$ (as in our case).

Assume that $X'$ and $X''$ are Cayley graphs  of the group~$A$, that are isomorphic to~$X$. Then $\aut(X')$ and $\aut(X'')$ are isomorphic to $\aut(X)\le \AGaL(1,q)$. By the above, any isomorphism from $X'$ to $X''$ induces an automorphism of $A$. Thus the number of graphs in Lemma~\ref{070825a} is at most $|\aut(A)|=|\GL(d,p)|\le p^{d^2}$.
\eprf

The order of the quantity  $\binom{n+\varepsilon}{m}$  is about
$\binom{\sqrt{q}}{\sqrt{q}/e}\cong p^{dp^{d/2}/{2e}}$, which for large $d$ is larger than $p^{d^2}$. By Lemmas~\ref{070825a} and~\ref{070825b}, this implies that there is a strongly regular graph $X'$ not isomorphic to $X$ and having the same parameters as~$X$. However, it is known that the $2$-dim WL algorithm does not distinguish any two strongly regular graphs with the same parameters. Thus the Weisfeiler-Leman dimension of the graph~$X$ is at least~$3$.

\bibliographystyle{amsplain}

\end{document}